\documentclass[11pt,twoside,reqno]{amsart}
\usepackage[all]{xy}
        \usepackage {amssymb,latexsym,amsthm,amsmath,mathtools,dsfont}
        \usepackage{enumitem,color}

\usepackage{mathtools}

\usepackage{mathrsfs}
\usepackage[
backend=biber,maxbibnames=999,doi=false,isbn=false,url=false
]{biblatex} 
\usepackage{hyperref}
\long\def\symbolfootnote[#1]#2{\begingroup%
\def\thefootnote{\fnsymbol{footnote}}\footnote[#1]{#2}\endgroup}

\newcommand{\F}{\ensuremath{\mathscr{F}}}

\newcommand{\GL}{\mathrm{GL}}
\newcommand{\SL}{\mathrm{SL}}
\newcommand{\SU}{\mathrm{SU}}
\newcommand{\Sp}{\mathrm{Sp}}
\newcommand{\SP}{\mathrm{Sp}}

\newcommand{\fq}{\ensuremath{\mathbb{F}_q}}

\newcommand{\U}{\mathrm{U}}

\makeatletter
\def\imod#1{\allowbreak\mkern10mu({\operator@font mod}\,\,#1)}
\makeatother
\makeatletter
\renewcommand*\env@matrix[1][*\c@MaxMatrixCols c]{%
  \hskip -\arraycolsep
  \let\@ifnextchar\new@ifnextchar
  \array{#1}}
\makeatother
\usepackage[capitalize,nameinlink]{cleveref} 
\usepackage{comment} 
\usepackage{xcolor} 
\hypersetup{
    colorlinks,
    linkcolor={red!80!black},
    citecolor={green!80!black},
    urlcolor={blue!80!black}
}
\makeatletter
\def\Ddots{\mathinner{\mkern1mu\raise\p@
\vbox{\kern7\p@\hbox{.}}\mkern2mu
\raise4\p@\hbox{.}\mkern2mu\raise7\p@\hbox{.}\mkern1mu}}
\makeatother

\newtheorem{theorem}{Theorem}[section]
\newtheorem{question}{Question}[section]
\newtheorem{lemma}[theorem]{Lemma}

\newtheorem*{theorem*}{Theorem}
\theoremstyle{definition}
\newtheorem{definition}[theorem]{Definition}

\newtheorem{example}[theorem]{Example}
\newtheorem{problem}[theorem]{Problem}
\numberwithin{equation}{section}
\newcommand{\ignore}[1]{}

\newcommand{\mynote}[1]{}
\begin{document}
\setcounter{section}{0}
\title{A survey on power maps in groups}
\author[Panja S.]{Saikat Panja}
\email{panjasaikat300@gmail.com}
\address{Stat-Math Unit, ISI Bangalore, 8th Mile, Mysore Rd, RVCE Post, Bengaluru, Karnataka 560059, India}
\author[Singh A.]{Anupam Singh}
\email{anupamk18@gmail.com}
\address{IISER Pune, Dr. Homi Bhabha Road, Pashan, Pune 411 008, India}
\thanks{The first named author (Panja) is supported by an NBHM fellowship for postdoctoral studies. The second-named author is funded by an NBHM research grant 2011/23/2023/NBHM(RP)/RDII/5955 for this research.}
\date{\today}


\begin{abstract}
The study of word maps on groups has been of deep interest in recent years. This survey focuses on the case of power maps on groups; \emph{viz.} the map $x\mapsto x^M$ for a group $G$, and an integer $M\geq 2$. Here, we accumulate various results on the subject and pose some questions.
\end{abstract}
\maketitle
\section{Introduction}\label{sec:intro}

Given a group $G$ and an element $g\in G$, the immediate operation which comes to one's mind is $g^2=g.g$, and if one keeps repeating the process $M-1$-times, we end up getting $g^M$. This defines a map $\theta_M\colon G\longrightarrow G$ defined as $\theta_M(g)=g^M$, and this map is known as \emph{power map}. This map is not a group homomorphism when $G$ is non-abelian. They appear naturally in several problems in mathematics, such as representation theory, Thompson's conjecture, discrete log problems, etc, as we see with different examples later on. 

Let us first describe a more general scenario called \emph{word map}. Given a word $w$ in the free group $\F_t$ on $t$ generators, and a group $G$, we get a map
\begin{align*}
    w\colon G^t\longrightarrow G,
\end{align*}
by plugging in the tuples from $G^t$, i.e., $(g_1, \ldots, g_t) \mapsto w(g_1, \ldots, g_t)$ which is called a word map. For example, when $w=x_1x_2x_1^{-1}x_2^{-1}\in \F_2$, which is known as commutator word, defines a map $G^2\longrightarrow G$, by $(g_1, g_2)\mapsto g_1g_2g_1^{-1}g_2^{-1}$. Ore conjectured in 1951 that for a finite non-abelian simple group, the commutator map is surjective after observing it for the alternating group in ~\cite{Ore51}. A related conjecture of Thompson (which implies Ore's conjecture) states that every finite non-cyclic simple group $G$ contains a conjugacy class $\mathscr{C}$ with $\mathscr{C}^2=G$, where $\mathscr{C}^2=\{xy \mid x,y\in \mathscr{C}\}$, see~\cite{Thompson61}. In the paper~\cite{Thompson61}, Thompson proved that Ore's conjecture holds for the groups $\mathrm{PSL}_n(q)$, the projective special linear groups of rank $n$, defined over a finite field of cardinality $q$.
Later by collective efforts in 2010, it was proved that Ore's conjecture holds true, see~\cite{LOST10} and the references therein. 

Hereafter, by a finite simple group, we will always mean~\emph{finite non-abelian simple group}. In the last 2-3 decades, the problems centred around the word maps have been an active area of research. The general theme pursued in this regard consists of the following questions:
\begin{question}\label{ques1}
Let $G$ be a group and $w$ be a word. Then
\begin{enumerate}
    \item Is it true that $w(G)=G$?
    \item If $w(G)\neq G$ is there a constant $k_w\in\mathbb{N}$ such that $w(G)^{k_w} = \langle w(G)\rangle$.
\end{enumerate}
\end{question}
These questions are also studied for a family of groups (e.g. finite simple groups, algebraic groups, Lie groups etc) and/or a family of words (e.g. $w$ being a power map or in $\F_1\backslash \F_2$ etc). Note that the image of a word map always (a) contains identity, and (b) is conjugacy-invariant (in fact, automorphism-invariant). Several interesting results are known when we restrict ourselves to finite non-abelian simple groups. In the year 2001, Liebeck and Shalev proved that for a finite simple group $G$ and an arbitrary word $w$, there exists a constant $c(w)$ such that $w(G)^{c(w)}=G$,~see \cite{LiebeckShalev01}. This was advanced to the following result of Shalev in 2009, where he proves that if $w\neq 1$ is a word, then there exists $N(w)$ such that for every finite simple group $G$ with $|G|> N(w)$ we have $w(G)^3=G$. Later, this was improved to the result $w(G)^2=G$ in \cite{LarsenShalevTiep11}. This is the best possible solution to the Waring problem for simple groups, as under power maps, the image need not be surjective. For example, the map $x\mapsto x^2$ is not surjective for the group $A_4$, the alternating group on $4$ letters, in fact, for any group having even order elements. Power maps also have the property of having dense image ratios, see \cite{panja2024ratios}.

Power maps have numerous applications, as we discuss below. 
Recall that a conjugacy class $\mathscr{C}$ of a group $G$ is called \emph{real} if $\mathscr{C}=\mathscr{C}^{-1}$.
The number of real conjugacy classes matches the number of real characters of $G$. Brauer's problem 14 asks whether the number of characters with the Frobenius-Schur indicator $1$ can be expressed in terms of group properties. 
A recent solution of this problem in~\cite{MurraySambale23} also shows that the number of real conjugacy classes is given by $s(2)/|G|$, where
$$s(2) = \left| \left\{(a,b)\mid a^2 b^2 = 1 \right\} \right|.$$
With some work, it can be shown that the number $s(2)$ can be found out by considering the \emph{fibers} of the word map $x\mapsto x^2$ and $x\mapsto x^4$, see~\cite{panja2024fibers}; where by fiber of a map $f\colon A\rightarrow B$ we mean the set $f^{-1}(b)$ for some $b\in B$. Furthermore, for groups whose irreducible characters are real, the column sums are given by the number of square roots of conjugacy class representatives, see~\cite{Isaacs76}.

A theorem of Martinez and Zelmanov \cite{MaZe96} from the 1990s shows that any element of a sufficiently large simple group can be written as a product of $f(k)$ many $k$-th powers. This was also independently proved by Saxl and Wilson in \cite{SaWi97}. This simple-looking result has far-reaching consequences. For example, in \cite{SaWi97}, the authors prove that if $G$ is a Cartesian product of nonabelian finite simple groups with $G$ being a 
finitely generated profinite group, then every subgroup $H$ of finite index in $G$ is open. For an element $g\in G$, the solutions of the equation $x^M=g$ have been studied in many contexts. For a $\mathbb{C}$-character of a group $G$ define the higher Frobenius-Schur indicator as
$$\varepsilon_m(\chi)=\dfrac{1}{|G|}\sum\limits_{g\in G}\chi(g^m).$$
Then it can be shown that for any $h\in G$, one has
$$\sum\limits_{\chi\in \mathrm{Irr}(G)}\varepsilon_m(\chi)\chi(h^{-1})=|\left\{y\in G|y^m=h\right\}|,$$
where $\mathrm{Irr}(G)$ denotes the set of irreducible characters of $G$, see \cite[Section 9]{Huppert98Book}. 

In this survey article, we bring together some of the results known about power maps. For finite classical groups, these are written in terms of generating functions, as is the case for the number of conjugacy classes for these groups. We also briefly mention the scenario for algebraic groups and Lie groups, including the more general contexts of word maps.

\section{Power maps in finite groups}\label{sec:finite}

One of the approaches to studying various problems for finite classical groups is the generating function method. If we have a family of groups $\{\mathcal G_n\}_{n\geq 1}$ and want to study the property $\mathcal P(\mathcal G_n)$, then the generating function for the same would be $\mathcal P_{\mathcal G}(z) = 1 + \displaystyle\sum_{n\geq 1} \mathcal P(\mathcal G_n)z^n$. Some examples of such a family of groups are the symmetric groups $S_n$, alternating groups $A_n$, and various classical groups: general linear groups, unitary, orthogonal and symplectic groups. Some examples of properties could be the number of conjugacy classes, number of regular semisimple elements, number of regular elements, number of semisimple elements, number of certain characters, number of squares, number of $M$-th powers, etc. A classic example is that for the symmetric group, the number of conjugacy classes is given by the partition function. Several probabilistic problems are formulated this way, and asymptotic results are obtained. Here, we will mostly talk about the power maps.  

Blum studied the number of square permutations in a symmetric group $S_n$ on $n$ letters in the paper \cite{Blum74}, where it is shown that if $f(n)$ denotes the number of square permutations in $S_n$, then the generating function 
\begin{align*}
    \sum\limits_{n=0}^{\infty} \dfrac{f(n)}{n!} z^n = \left[\dfrac{1+z}{1-z}\right]^{1/2} \prod\limits_{k=1}^{\infty}\cosh{\left(\dfrac{z^{2k}}{2k}\right)}
\end{align*}
It was further shown that $f(2k+1)=(2k+1)f(2k)$. 
Later in \cite{MikAndDen00}, the authors proved that the probability $p_r(n)$ that a random permutation from $S_n$ has an $r$-th root, with $r$ prime, is monotonically nonincreasing in $n$. Let $p_n(m)$ denote the proportion of symmetric group $S_n$ permutations that admit an $m$-th root. Then it has been shown that
\begin{align*}
  \lim_{n\to \infty}  p_n(m)\,\,\,\, \sim \,\,\,\,\dfrac{\pi_m}{n^{1-\varphi(m)/m}},
\end{align*}
for an explicit constant $\pi_m$ and $\varphi$ being the Euler function, see \cite[Section 4.3, Theorem]{Pouyanne02}. There has also been study of finding the number of $M$-th root in these cases, for example, \cite{LeanosMorenoRivera12}, \cite{GlebskyLiconRivera23}.
Results of similar flavours for the wreath product of groups appear in \cite{KunduMondal22}. The methods follow a similar pattern for all of the above cases, as we discuss below.

Note that an element $g\in G$ is an $M$-th power if and only if for any $h\in G$, the element $hgh^{-1}$ is also an $M$-th power. Thus, to identify elements that are $M$-th powers, it is sensible to work with the conjugacy classes. Hence, a way to enumerate $M$-th power elements in a finite group might consist of the following steps:
\begin{enumerate}
\item identify the conjugacy classes which are $M$-th power,
\item calculate the size of the centralizers corresponding to each class (\emph{or} equivalently one should know the size of the conjugacy classes)
\item find an efficient way to combine this information if one encounters a family of groups (say, the family of finite general linear groups $\GL(n, q)$).
\end{enumerate}
In a series of three articles~\cite{KunduSingh24, PanjaSingh22symporth, PanjaSinghUnitary2024} $M$-th power maps have been described in the case of classical groups. We briefly describe these results below.

In what follows, $\Phi$ will denote the set of all monic irreducible polynomials $f(t)$ over $\fq$, presumably $t\not\in\Phi$, and $\Lambda$ will denote the set of all partitions of all $n\in\mathbb{N}$. The \textit{dual} of a monic degree $r$ polynomial $f(t)\in k[t]$ satisfying $f(0)\neq 0$, is the polynomial given by $f^*(t)= f(0)^{-1}t^rf(t^{-1})$. A \emph{self reciprocal irreducible monic} (SRIM) polynomial is a monic irreducible polynomial satisfying $f=f^*$. A \emph{self conjugate irreducible monic} (SCIM) polynomial is a polynomial $f$ of degree $d$ such that $\widetilde{f}=f$ where $\widetilde{f}(t) = \overline{f(0)}^{-1} t^d \bar f(t^{-1})$. Since the conjugacy classes of classical groups are described in terms of certain polynomials and partitions, the same can be done for the powers.
\begin{lemma}\label{lem:root-single-block} The $M$-th root in the following matrix groups are given as follows: 
\begin{enumerate}
\item \cite[Proposition 4.5]{KunduSingh24} Let $A\in\GL(n, q)$ has characteristic polynomial $f$ of degree $n$, which is irreducible. Then, $\alpha^M = A$ has a solution in $\GL(n, q)$ if and only if $f(t^M)$ has an irreducible factor of degree $n$.
\item \cite[Lemma 5.1 and Lemma 5.2]{PanjaSingh22symporth} Let $A\in\Sp(2n,q)$ has characteristic polynomial $f$, which is {\color{red}SRIM} of degree $2n$. Then $\alpha^M = A$, has a solution in $\Sp(2n,q)$, if and only if $f(t^M)$ has a SRIM factor of degree $2n$.
\item \cite[Lemma 4.1 and Lemma 4.2]{PanjaSinghUnitary2024}Let $A\in \U(n,q^2)$ has characteristic polynomial $f$, which is a {\color{red}SCIM} polynomial. Then $\alpha^M=A$ has a solution in $\U(n, q^2)$  if and only if $f(t^M)$ has a SCIM factor of degree $n$.
\end{enumerate}
\end{lemma}
This motivates us to study several classes of polynomials, as we have noted down here.
\begin{definition}
Let $f(t)$ be an irreducible polynomial over $\mathbb{F}_q$. Then,
\begin{itemize}
\item the polynomial $f$ is said to be an $M$-power polynomial if $f(t^M)$ has an irreducible factor of degree $\deg f$;
\item if $f$ is a self-reciprocal polynomial then $f$ is said to be an $M^*$-power polynomial if $f(t^M)$ has a SRIM factor of degree $\deg f$; and
\item  if $f$ is a self-conjugate polynomial then $f$ is said to be an $\widetilde{M}$-power polynomial if $f(t^M)$ has a SCIM factor of degree $\deg f$.
\end{itemize}
\end{definition}
It may seem that the existence of a root in $\GL(m, q)$ is the same as the existence of roots in other groups, but this is not the case, as the following examples illustrate. Using the previous lemma, we give examples in terms of the appropriate polynomials, as that suffices to prove our point. 
\begin{example}
Consider $\mathbb{F}_5$ and the polynomial $t^4 + 3t^3 + t^2 + 3t + 1 \in\mathbb{F}_5[t]$. Then $$f(t^2)=t^8 + 3t^6 + t^4 + 3t^2 + 1 = (t^4 + 2t^3 + t^2 + 3t + 1)(t^4 + 3t^3 + t^2 + 2t + 1).$$ 
Thus, it is a $2$-power polynomial but not a $2^*$-power SRIM polynomial.
\end{example}
\begin{example}
Let $\mathbb{F}_{25}=\mathbb{F}_5[a]$ and consider the polynomial $f(t)=t^5 + a^5t^4 + t^3 + a^2  t^2 + at + 1$. Then 
\begin{align*}
f(t^3) &= (t^2 + t + 3a )\cdot (t^2 + (2a + 1)t + 3a + 2 )\cdot (t^2 + (3a + 3)t + 4a + 3) \\
&\cdot (t^3 + t^2 + (3a + 3)t + 4a + 1) \cdot (t^3 + (2a + 1)t^2 + (2a + 1)t + 4a + 1) \\&\cdot (t^3 + (3a + 3)t^2 + t + 4a + 1), \end{align*}
and none of the irreducible factors of degree $3$ are self-conjugate. Thus, $f$ is a $3$-power polynomial but not a $\widetilde{3}$-power polynomial.
\end{example}
The number of such polynomials can be found using combinatorial techniques. It is well known that for a SRIM polynomial, the degree should be even, whereas the degree of a SCIM polynomial is always odd. We note down some of these here. 
\begin{lemma} Let $\mu$ denote the M\"{o}bius function. Then one has the following:
\begin{itemize}
\item\cite[Proposition 3.3]{KunduSingh24} Let $N_M(q, d)$ denotes the number of $M$-power polynomials of degree $d$. Then for $d>1$
 $$N_M(q,d)=\dfrac{1}{d}\sum\limits_{l \mid d}\mu(l)\dfrac{(M(q^{d/l}-1), q^d-1)}{(M,q^d-1)}.$$
\item\cite[Proposition 4.6]{PanjaSingh22symporth} Let $N^*_M(q, 2d)$ denotes the number of $M^*$-power SRIM polynomial of degree $2d$, $d\geq 1$. Then 
$$N^*_M(q,2d)=\dfrac{1}{2d}    \sum\limits_{\substack{l=\textup{odd}\\l\mid 2d}}\mu(l)\dfrac{(M(q^{2d/l}-1),q^d+1)}{(M,q^{2d}-1)}.$$
\item \cite[Lemma 3.2]{PanjaSinghUnitary2024}
For odd $d$, let $\widetilde{N}_{M}(q,d)$ denote the number of $\widetilde{M}$-power SCIM polynomials of degree $d$. Then we have 
\begin{align*}
\widetilde{N}_{M}(q,d) = \dfrac{1}{d}\sum\limits_{l\mid d}^{}\mu(l)\dfrac{\left(M(q^{2d/l}-1),q^d+1\right)}{(M,q^{2d}-1)}.
\end{align*}
\end{itemize}   
\end{lemma}
\begin{example}
Let us take $M=2$ and $q$ to be odd. In this case, we have the following values;
\begin{table}[h]
\centering
\begin{tabular}{|c|c|c|c|} \hline
$d$ & $N_2(q,d)$ & $N^*_2(q,2d)$ & $\widetilde{N}_2(q,2d-1)$ \\    \hline
$2$ & $\dfrac{1}{2}\left(q^2-q\right)$ & $\dfrac{1}{8}\left(q^2+1\right)$ & $\dfrac{1}{6} \left( q^3 - q \right) $\\
$3$ & $\dfrac{1}{2}\left(q^3-q\right)$ & $\dfrac{1}{12}\left(q^3-q\right)$ & $\dfrac{1}{10} \left(q^5 + 1 \right)$\\
$6$ & $\dfrac{1}{2}\left(q^6-q^3-q^2+q\right)$ & $\dfrac{1}{24} \left( q^6 + q^2 - 2 \right)
$ & $\dfrac{1}{22} \left(  q^{11} - 1 \right)$ \\ \hline \end{tabular} \caption{Some values of $N_M$, $N^*_M$ and $\widetilde{N}_M$.}
    \label{tab:values-of-N-M}
\end{table}
\end{example}
Recall that a matrix is called \emph{separable} if all its eigenvalues are distinct. Thus, using \cref{lem:root-single-block}, the generating functions for the proportion of regular semisimple elements which admit an $M$-th root can be derived. We note them down here.
\begin{theorem}\label{thm:regular-semisimple}
\begin{enumerate}
\item\cite[Theorem 5.3]{KunduSingh24} Let $S_\GL^M(n,q)$ denote the probability of an element of $\GL(n, q)$ being an $M$-th power separable element. Let $S_\GL(z)=1+\sum\limits_{n=1}^\infty S_\GL^M(n,q)z^n$. Then we have, 
\begin{align*}
S_\GL(z)=\prod_{d\geq 1}\left(1+\dfrac{z^d}{q^d-1}\right)^{{N}_M(q,d)}.
\end{align*}
\item\cite[Theorem 5.6]{PanjaSingh22symporth} Let $S^M_{\Sp}(n, q)$ denote the probability of an element to be $M$-power separable in $\Sp(2n,q)$ and 
$S_{\SP}(z) = 1+\sum\limits_{m=1}^{\infty}S^M_{\SP}(m,q)z^m$. Then 
\begin{equation*}
S_{\SP}(z) = \displaystyle\prod_{d=1}^{\infty}\left(1+\dfrac{z^d}{q^d+1}\right)^{N_M^*(q,2d)}\prod_{d=1}^{\infty}\left(1+\dfrac{z^d}{q^d-1}\right)^{R_M^*(q,2d)}.
\end{equation*}
\item\cite[Theorem 4.4]{PanjaSinghUnitary2024} Let $S_\U^M(n, q)$ denote the probability of an element of $\U(n, q^2)$ being an $M$-th power separable element. Let $S_\U(z) = 1+ \sum\limits_{n=1}^\infty S_\U^M(n,q)z^n$. Then 
\begin{align*}
S_\U(z)=\prod_{\substack{d\text{ odd}\\d\geq 1}}\left(1+\dfrac{z^d}{q^d+1}\right)^{\widetilde{N}_M(q,d)}\prod_{d\geq 1}\left(1+\dfrac{z^{2d}}{q^{2d}-1}\right)^{\widetilde{R}_M(q,d)}.
\end{align*}
\end{enumerate}
Here $R^*_M(q, 2n)$ denotes the number of pairs $\{\phi,\phi^*\}$, where $\phi~(\neq\phi^*)$ is an irreducible monic polynomial of degree $n\geq 2$ and $\phi$ is an $M$-power polynomial; which means \begin{align*}R^*_M(q,2n)=\begin{cases}
\frac{1}{2}N_M(q,n) & n\text{ is odd}\\
\frac{1}{2}\left(N_M(q,n)-N_M^*(q,n)\right)& n\text{ is even}
\end{cases}.
\end{align*}
Also, by definition $$\widetilde{R}_M(q,d)=\dfrac{N_M(q^2,d)-\widetilde{N}_{M}(q,d)}{2},$$ to be the number of pairs of $g,~ \widetilde{g}\in\Phi_M$ such that $g\neq \widetilde{g}$, where $\widetilde{N}_M(q,d)=0$ if $d$ is even.
\end{theorem}
There are results which provide generating functions for the proportion of various other elements admitting an $M$-th root of the given groups. We refer an interested reader to the articles \cite{KunduSingh24}, \cite{PanjaSingh22symporth}, and \cite{PanjaSinghUnitary2024} for further details.

Now that we are able to identify some class of elements admitting an $M$-th root, it is natural to ask, given an element $g\in G$ having an $M$-th root, how many $M$-th roots do exist for $g$? This question has been answered for the identity elements in finite classical groups. In the article \cite{panja2024roots}, we derive the generating functions for the probability of an element satisfying $x^M=1$, where $M\geq 2$ is an integer. We call these elements the \emph{$M$-th roots of identity}. Note that if an element is a root of identity, then all its conjugates are as well. Thus, we first find the conjugacy class representatives, which are $M$-th root of identity, and proceed similarly as above. For example, we have the following theorem;
\begin{theorem}\label{th:GL-all-probability}\cite[Theorem 4.3]{panja2024roots}
Let $a_n$ denote the number of elements in $\GL_n(q)$ which are $M$-th root of identity. Let $M=t\cdot p^r$, where $p\nmid t$. Then the generating function of the probability $a_n/|\GL_n(q)|$ is given by
\begin{align*}
&1+\sum\limits_{n=1}^{\infty}\dfrac{a_n}{|\GL_n(q)|}z^n=
\prod\limits_{d|t}\left(1+\sum\limits_{m\geq 1}\sum\limits_{\substack{\lambda\vdash m\\\lambda_1\leq p^r}} \dfrac{z^{me(d)}}{q^{e(d)\cdot(\sum_{i}(\lambda_i')^2)}\prod\limits_{i\geq 1}\left(\dfrac{1}{q^{e(d)}}\right)_{m_i(\lambda_{\varphi})}}\right)^{\frac{\phi(d)}{e(d)}},
\end{align*}
    where $e(d)$ denotes the multiplicative order of $q$ in $\mathbb Z/d\mathbb Z^\times$. Here the symbol 
$\left(\dfrac{u}{q}\right)_i$ denotes the quantity $\left(1-\dfrac{u}{q}\right)\left(1-\dfrac{u}{q^2}\right)\ldots \left(1-\dfrac{u}{q^i}\right)$.
\end{theorem}
\begin{example} This is snippet from \cite[Section 5]{panja2024roots}. 
Let $M\neq 2$ be a prime and $q\equiv-1 \pmod{M}$, for example, $(q, M)=(41,7)$. Let $b_n$ denote the proportion of $M$-th roots of identity in $\GL(n, q)$. Then using we get 
\begin{align*}
1+\sum\limits_{n=1}^{\infty}b_nz^n&=\left(1+\sum\limits_{m=1}^{\infty}\dfrac{z^m}{|\GL_m(q)|}\right)\left(1+\sum\limits_{m=1}^{\infty}\dfrac{z^{2m}}{|\GL_m(q^2)|}\right)^{{(M-1)}/{2}}.
\end{align*}
We now divide the computation into two cases. We first consider odd values of $n$. In this case, we should have an odd power of $z$, coming from the first term of the product. Other contributing powers of $z$ will have all even power. Hence, the probability of being an $M$-th root is
{
\begin{align*}
\sum\limits_{\substack{1\leq j \leq M\\j=\text{odd}}}\dfrac{1}{|\GL_j(q)|}\cdot\left(\sum\limits_{\lambda\vdash \frac{M-j}{2}}\prod\limits_{\ell}\dfrac{1}{|\GL_{\lambda_{\ell}}(q^2)|}\right),
\end{align*}
where $\ell$ runs over the subscripts of the parts of $\lambda=(\lambda_1,\lambda_2,\ldots)$.
When $n$ is even, using the same argument as before, we get the resulting probability to be
\begin{align*}
\sum\limits_{\substack{0\leq j \leq M\\j=\text{even}}}\dfrac{1}{|\GL_j(q)|}\cdot\left(\sum\limits_{\lambda\vdash \frac{M-j}{2}}\prod\limits_{\ell}\dfrac{1}{|\GL_{\lambda_{\ell}}(q^2)|}\right),
\end{align*}
where $\ell$ runs over the subscripts of the parts of $\lambda=(\lambda_1,\lambda_2,\ldots)$ and $|\GL_0(q)|$ is $1$ by convention.}
\end{example}
Towards the end of this section, based on the results discussed above, we propose the following problems to look into.
\begin{problem}
Find the generating functions for the exact number of elements admitting $M$-th root in all classical groups.
\end{problem}
\begin{problem}
Given a finite exceptional group of Lie type, classify the elements admitting an $M$-th root.
\end{problem}
\begin{problem}
In the memoir \cite{FulmanNeumanPraeger05}, there are asymptotic values of proportions of certain classes of elements are given. For example, \cite[Theorem 2.1.4]{FulmanNeumanPraeger05} states that the limiting probability of proportion of separable elements in $\U(n,q^2)$ as $n\longrightarrow\infty$, denoted by $S_\U(\infty,q)$ satisfies
\begin{align*}
1-\dfrac{1}{q}-\dfrac{2}{q^3}+\dfrac{2}{q^4}<S_\U(\infty,q)
 <1-\dfrac{1}{q}-\dfrac{2}{q^3}+\dfrac{6}{q^4}.
\end{align*}
This motivates us to ask for limiting probability for the case of $M$-th power elements as well.
\end{problem}

In the article \cite{KulKundSingh22}, the authors investigate the asymptotic proportion of powers of regular semisimple, semisimple, and regular elements of a finite reductive group when $q\longrightarrow\infty$. For example, if $G(\mathbb{F}_q)^M$ denotes the image of the map $\theta_M\colon G(\mathbb{F}_q)\longrightarrow G(\mathbb{F}_q)$ given by $x\mapsto x^M$, then one has
\begin{align*}
\lim\limits_{q\longrightarrow\infty} \dfrac{|G(\mathbb{F}_q)^M|}{|G(\mathbb{F}_q)|} = \sum\limits_{{T}(\mathbb{F}_q)=T_{d_1,d_2,\ldots,d_k}}\dfrac{1}{\mathbb{W}_{T(\mathbb{F}_q)}(M,d_1)\ldots(M,d_k)},
\end{align*}
where $T(\mathbb{F}_q) = C_{d_1}\times\ldots \times C_{d_k}$ and $W_{{T}(\mathbb{F}_q)}$ denotes the Weyl group of the corresponding torus.
\section{Power maps: Beyond finite groups}\label{sec:algebraic}

 Let $k$ be an algebraically closed field and $G$ be an algebraic group over $k$. For any non-trivial word $w\in \F_t$, we get an algebraic morphism $w$ on $G$. In \cite{Bo83} (see Theorem B and Theorem 1), Borel (what is now called Borel dominance theorem named after him) showed that when $G$ is a connected semisimple algebraic group over a field $k$, the image $w(G)$ is dominant. In view of this, one can show that $w(G)^2 = G$ for such groups. Thus, for the power map, the square of the image would be the whole $G$. Further extension of this result is studied in~\cite{GKP16}, \cite{GKP18}. However, for particular words like power maps, better results are known.

 Consider the power map $\theta_M \colon G \rightarrow G$, an algebraic morphism. In 2002-03, Chatterjee~\cite{Chatterjee02} and~\cite{Chatterjee03} proved that for connected semisimple algebraic groups, $\theta_M$ is mostly surjective. This was also shown by Steinberg~\cite{Steinberg03} following an independent method. More precisely, they proved the following:
 \begin{theorem}
 Let $G$ be a connected semisimple algebraic group over $k$. Then, the map $\theta_M$ is surjective if and only if  $M$ is coprime to $bz$ where $b$ is a bad prime and $z$ is the size of the center.  
 \end{theorem}
 The bad primes, in the case of $G$ simple algebraic group, are as follows:
 \vskip5mm
 \begin{center}
 \begin{tabular}{|c|c|c|} \hline
 Group type & Bad primes & Order of center\\  \hline
  $A_l, l\geq 1$    & $1$ & $l+1$ \\
 $B_l, l\geq 2; C_l, l\geq 3; D_l, l\geq 4$     & $2$ & $2$\\ 
 $G_2, F_4$ &$2, 3$ & $1$\\
 $E_6$ & $2,3$ & $3$\\
 $E_7$ & $2,3$ & $2$\\
 $E_8$ & $2,3, 5$ & $1$\\ \hline 
 \end{tabular}
 \end{center}
 \vskip2mm

 \noindent One could use this data and immediately figure out that over an algebraically closed field, the power maps are surjective when the characteristic of the field is large enough. However, the story gets more complicated when either field is not algebraically closed or the group is not semisimple.  
 \begin{example}
 The square map $\theta_2$ is surjective on $\GL_2(\mathbb C)$ but it is not surjective on $\GL_2(\mathbb R)$ and $\GL_2(q)$. 
 \end{example} 
 The surjectivity of the power map for $p$-adic algebraic groups and real algebraic groups is studied in \cite{Chatterjee09, Chatterjee11} and its relation with exponentiality problem in Lie groups is explored too. The density of images for power maps in certain solvable groups and Lie groups is explored in \cite{BM18, DM17}. Hui, Larsen and Shalev in \cite{HLS15} and Egorchenkova and Gordeev in \cite{EG19} have considered the generalization of Borel's result over the base field for any word.

 Avni et al. in \cite{AGKS13} have studied word map questions for simple algebraic groups over the $p$-adic integers and general local rings. For example, they prove that if $G$ is a semisimple, simply connected algebraic group over $\mathbb Q$, and $w$ is a nontrivial word, then for large enough primes $p$, $w(G(\mathbb Z_p))^3=G(\mathbb Z_p)$. 
In \cite{AM19} Avni and Meiri, look at the word map for $\SL_n(\mathbb Z)$. Among several interesting results, they show that (Theorem 1.1) there exists a constant $C$ (which can be taken to be $87$) so that, for any word $w$, there is a positive integer $d(w)$ such that for all $d > d(w)$, every element of $\SL_d(\mathbb Z)$ is a product of at most $C$ elements of the image $w(\SL_d(\mathbb Z))$.

In \cite{ET14}, Elkasapy and Thom looked at the word map in two variables on the group $\SU(n)$ and showed that for any word not in $[\F_2,\F_2]$, the word map is surjective for large enough $n$.
In \cite{BZ16, GG20, KM22,JS21}, some more words are looked at on the groups $\rm{PGL}_2$, $\rm{PSL}_2$, $\SL_2(\mathbb C)$ and $\SU(2)$.

\printbibliography
\end{document}